\documentclass[10pt]{amsart}

\usepackage{preamble/Pakete}
\usepackage{preamble/Befehle}

\begin{document}

\title[Koszul duality and determinant formulas]{Koszul duality and
determinant formulas for affine hyperplane arrangements}
\author{Jens Niklas Eberhardt}
\author{Carl Mautner}
\date{9 November 2026}

\begin{abstract}
For real affine hyperplane arrangements, we prove Lam's determinant
conjecture for the inverse Varchenko matrix restricted to bounded chambers.  The proof combines Lam's comparison
of twisted intersection forms with the period determinant formula of
Varchenko, Loeser--Sabbah and Douai--Terao.
For normal crossing arrangements, we introduce the taut path algebra and the chamber intersection algebra and establish a Koszul duality between
them, yielding an
explicit multivariate formula for the inverse Varchenko matrix.
After specializing the parameters to $q$, its restriction to bounded chambers
is a scalar multiple of the $q$-intersection matrix, proving a determinant
conjecture of the authors.  We also discuss how the Koszul duality
results generalize to affine oriented matroids.
\end{abstract}
\maketitle

\section{Introduction}
\subsection{Two matrices attached to a hyperplane arrangement}
Let $\C=\{H_i\}_{i\in E_0}$ be an essential arrangement of distinct
affine hyperplanes in $\R^r$. We denote its chambers by $\cF$ and its
bounded chambers by $\cB\subseteq\cF$.

For a finite polytopal complex $F$, the $h$-polynomial of $F$ is defined as
\[
 h(F,t)=\sum_{\varnothing\ne G\leq F}(t-1)^{\dim G},
 \qquad h(\varnothing,t)=0,
\]
where $G$ runs over its nonempty faces. For a closed face $F$, let
$F_\bd$ be its complex of bounded faces. For chambers
$P,Q\in\cF$, let $\Sep(P,Q)$ be the set of hyperplanes separating
them and put $d(P,Q)=|\Sep(P,Q)|$. 

Consider the polynomial rings $\ZM[q]$ and $\ZM[\mathbf b] = \ZM[b_i~|~i\in E_0]$ in variables $q$ and $b_i$ for $i \in E_0$ respectively.   As shorthand, for a subset $J \subseteq E_0$, we write $b_J=\prod_{i\in J}b_i$.

Our main results involve two $\cF\times \cF$-matrices. 

\begin{definition}\cite{Varchenko1993} The \textit{Varchenko matrix} $V_b$ is the $\cF\times \cF$-matrix with $(P,Q)$-entry
\[
 V_b(P,Q)=b_{\Sep(P,Q)} \in \ZM[\mathbf b],
\]
for any pair of chambers $P,Q\in\cF$.  Let $V_q=V_b|_{b_i=q}$ be the matrix with entries in $\ZM[q]$ obtained from $V_b$ by specializing all of the variables $b_i$ to $q$.
\end{definition}

\begin{definition} The \textit{$q$-intersection matrix} $W_q$ is the $\cF\times \cF$-matrix with $(P,Q)$-entry
\[
W_q(P,Q)=(-q)^{d(P,Q)}h((\bar P\cap\bar Q)_\bd,q^2) \in \ZM[q],
\]
for any pair of chambers $P,Q\in\cF$.  
\end{definition}

For an $\cF \times \cF$-matrix $X$, write $X^{\cB}$ for its
$\cB\times\cB$-submatrix. For hyperplane arrangements with normal crossings, the $\cB \times \cB$-matrix $S_q:=W_q^{\cB}$ is the $q$-intersection matrix introduced previously in \cite{EM}.

Our main theorems are a determinant formula for $(V_b^{-1})^{\cB}$, the $\cB \times \cB$-submatrix of the inverse of the Varchenko matrix, confirming a conjecture of Lam \cite[Conjecture 8.10]{Lam}, and a determinant formula for $S_q$ for normal crossing arrangements, confirming our conjecture~\cite[Conjecture 1.7]{EM}.

\subsection{Determinant Formulas}
Denote by $\widehat{\C}$ and $\C_0$ the homogenized and central
arrangements attached to $\C$, and by $M$ and $M_0$ the matroids of
their normal vectors. Their ground sets are $E=E_0\sqcup\{0\}$ and
$E_0$, respectively, where $0$ labels the hyperplane at infinity.
Write $\cL = \cL(M)$ for the lattice of flats of $M$, $\cL_0 \subseteq \cL$ for the set of flats containing $0$,
and $\Z(M_0)$ for the cyclic flats of $M_0$.

For a matroid $N$ we write $\beta(N)$ for Crapo's beta invariant and $\mup(N)$ for the unsigned Möbius number.  For any flat $F$ of $N$, we set \[w_F=\beta(N^F)\beta(N_F),\] where $N^F=N|F$ denotes the restriction to $F$ and $N_F=N/F$ denotes the contraction by $F$.  We recall additional matroid notation in Section~\ref{sec:background}.

The following theorem proves Lam's determinant conjecture~\cite[Conjecture 8.10]{Lam} in the case of real affine hyperplane
arrangements.
\begin{theorem}\label{thm:Lam-real}
For any essential real affine hyperplane arrangement,
\begin{equation}\label{eq:Lam-determinant}
 \det\bigl((V_b^{-1})^{\cB}\bigr)
 =\frac{\displaystyle\prod_{F\in\cL_0\setminus\{E\}}
          (1-b_{E\setminus F}^2)^{w_F}}
        {\displaystyle\prod_{F\in\cL\setminus(\cL_0\cup\{\varnothing\})}
          (1-b_F^2)^{w_F}}.
\end{equation}
\end{theorem}

Lam states his conjecture in terms of the bounded block of the Betti homology matrix, which he proves is equal to $(-1)^r(V_b^{-1})^{\cB}$, as recalled in
Section~\ref{subsec:betti-pairing}.

We also prove the determinant formula conjectured in
\cite[Conjecture 1.7]{EM}.  To state our result, we write
\[
 \qint{n}=1+q^2+\cdots+q^{2n-2}=\frac{1-q^{2n}}{1-q^2}.
\]

\begin{theorem}\label{thm:intersection-real}
If $\C$ is an essential real affine hyperplane arrangement with normal crossings, that is, if any non-empty intersection of $s$ hyperplanes in $\C$ has codimension $s$, then
\begin{equation}\label{eq:Sq-determinant}
 \det S_q
 =\prod_{\substack{K\in\Z(M_0)\\K\ne E_0}}
   \qint{|E_0\setminus K|}^{\beta(M_{0,K})\mup((M_0^K)^*)}.
\end{equation}
\end{theorem}

\begin{remark}
Specializing $q$ to 1, the entries of $S_1=S_q|_{q=1}$ are, up to sign, counts of vertices of common
faces. The $q=1$ determinant formula was proved in \cite[Theorem 1.4]{EM}
by identifying $S_1$ with the Brylawski--Schechtman--Varchenko form~\cite{SchechtmanVarchenko1991,BrylawskiVarchenko1997}
for the dual matroid $M_0^*$, with all weights equal to one.
For a unimodular arrangement, the bounded chambers index the toric
components of the compact core of a smooth hypertoric variety.
The matrix $S_1$ describes the intersection form on the top cohomology of this core
\cite{HauselSwartz2006}.

For realizable (and more generally orientable) matroids, the matrix $S_1$ can be used to study the representation theory of the matroidal Schur algebras of~\cite{Mautner,BradenMautner}.  We expect the existence of matroidal $q$-Schur algebras for which $S_q$ should play an analogous role.
\end{remark}

The proof of Theorem~\ref{thm:Lam-real} depends on a determinant formula for a de Rham pairing and many other results in the monumental work of Lam \cite{Lam} as well as the beautiful formulas for the determinants of period matrices due to Varchenko, Loeser--Sabbah and Douai--Terao \cite{Varchenko1990,LoeserSabbah1991,DouaiTerao1997}. We deduce Theorem~\ref{thm:intersection-real} by relating, for arrangements with normal crossings, the $q$-intersection matrix $W_q$ to the inverse of the Varchenko matrix, as we describe below.

\begin{remark}
Both determinant formulas are conjectured to hold in the more general setting of affine oriented
matroids, with an analogous normal crossing hypothesis for the second (see Definition~\ref{def:normal-crossing}). Our proofs
in the setting of hyperplane arrangements use realizability through the period determinant formula.
\end{remark}

\subsection{Koszul duality and the inverse Varchenko matrix}  To show that Theorem~\ref{thm:Lam-real} implies Theorem~\ref{thm:intersection-real}, we prove that when $\C$ has normal crossings, then the $q$-Varchenko matrix $V_q$ and the $q$-intersection matrix $W_q$
are inverse up to the scalar factor $(1-q^2)^r$.
In fact, we show the following stronger, categorical statement. Namely, we define the \emph{taut path algebra}
$\Taut$ and the \emph{intersection algebra} $\Inter$, and show that their Hilbert
matrices are $V_q$ and $(1-q^2)^{-r}W_{-q}$, respectively.
These algebras are Koszul dual by Theorem~\ref{thm:chamber-koszul}.
The Hilbert-series identity for Koszul dual algebras then gives the
following comparison.
\begin{theorem}\label{thm:intro-inverse}
For a normal crossing arrangement $\C$, we have
\[
 W_q=(1-q^2)^rV_q^{-1},\qquad
 S_q=(1-q^2)^r(V_q^{-1})^{\cB}.
\]
\end{theorem}
We prove this comparison more generally for normal crossing affine oriented matroids (see Definition~\ref{def:normal-crossing}).
In this setting as well, it reduces the conjectural determinant formula
\eqref{eq:Sq-determinant} to \eqref{eq:Lam-determinant} in the normal
crossing case.

\begin{remark}
The algebras $\Taut$ and $\Inter$ are closely related to algebras introduced by Braden--Licata--Proudfoot--Webster
\cite{BLPWgale,BradenLicataProudfootWebster2012}. In the realizable
normal crossing case, $\Inter$ appears in~\cite[Section 4.5]{BLPWgale} as the extended equivariant algebra
$\widetilde B_{\mathrm{ext}}$ defined on all feasible chambers.  The Koszul duality of Braden--Licata--Proudfoot--Webster concerns category $\mathcal O$ for Gale dual polarized arrangements.
Our theorem gives a duality on the full chamber set of a single affine
oriented matroid. For real arrangements, $S_q$ is also the graded
Euler form on the simple objects indexed by bounded chambers in
hypertoric category $\mathcal O$ \cite[Theorem 3.1]{EM}.
\end{remark}

\subsection{Outlook}\label{subsec:outlook}

\begin{enumerate*}[label=(\arabic*), mode=unboxed, itemjoin=\endgraf\noindent]
\item
It would be interesting to extend the Koszul duality to hyperplane arrangements (and affine oriented matroids) that do not have normal
crossings. The natural replacement for the intersection algebra is
\[
 \operatorname{Ext}_{\Taut}^{\bullet}(\Taut_0,\Taut_0),
\]
together with its internal multigrading and higher multiplications.
Taking alternating sums in cohomological degree recovers $V_b^{-1}$,
thereby categorifying Lam's Betti homology matrix at $b_0=0$, up to
the sign $(-1)^r$. For the arrangements describing $M_{0,n}$, the
bounded-chamber entries are string-corrected biadjoint scalar amplitudes
up to normalization \cite[Section 19.4]{Lam}, so this also gives a
categorification of these amplitudes. An explicit combinatorial
description and an extension to general $b_0$ would be desirable.

\item
A combinatorial proof avoiding the period determinant formula of
Loeser--Sabbah and Douai--Terao could establish both determinant formulas for
affine oriented matroids, with the normal crossing assumption for
the $q$-intersection formula. We already verified the latter formula
for the non-realizable pseudoline arrangement pictured on the cover
of \emph{Oriented Matroids}
\cite{OrientedMatroids1999}, and for a generic affine lifting of an
orientation of the V\'amos matroid. See \cite[Sections A.2--A.3]{EM}.

\item
It would also be interesting to apply our methods to the Koszulity
question raised by Kowalenko--Mautner for the algebras associated to
generic oriented matroid programs \cite[Section 1.3]{KowalenkoMautner}.
They establish quadratic duality in general and Koszulity for
Euclidean programs. An affirmative answer would extend the Koszul
duality of Braden--Licata--Proudfoot--Webster \cite{BLPWgale} to arbitrary
generic oriented matroid programs.
\end{enumerate*}

\subsection{Structure}
Section~\ref{sec:background} fixes notation for arrangements and affine
oriented matroids. Section~\ref{sec:periods} recalls the twisted
intersection forms and period determinant formula used in
Section~\ref{sec:determinants} to prove Theorem~\ref{thm:Lam-real}.
Section~\ref{sec:comb} collects the combinatorial lemmas for normal
crossing affine oriented matroids. Section~\ref{sec:koszul} proves
Koszul duality and the inverse-matrix identity.
Section~\ref{sec:qdet} deduces Theorem~\ref{thm:intersection-real}.

\section{Arrangements and affine oriented matroids}\label{sec:background}

\subsection{Affine and central arrangements}
Write $\C=\{H_i\}_{i\in E_0}$, where
$H_i=\{f_i=0\}$ and $f_i(x)=\ell_i(x)+c_i$ for linear functionals $\ell_i \in (\R^r)^*$ and scalars $c_i \in \R$. To $\C$ we attach the \emph{central arrangement} $\C_0=\{H_i^0\}_{i\in E_0}$ where $H_i^0=\{\ell_i=0\}$. Moreover, by homogenizing the defining equations of $\C$ we obtain the central arrangement $\widehat\C$ in $\R^{r+1}$ with hyperplanes defined by $\widehat f_i(x,t)=\ell_i(x)+c_it$ for $i\in E_0$ and an additional hyperplane at infinity defined by $\widehat f_0(x,t)=t$ and indexed by $0$. The slice $t=1$
recovers $\C$, and restriction of the forms $\widehat f_i$, $i\in E_0$,
to $t=0$ gives $\C_0$.

To a central arrangement, one may associate a matroid defined on the set indexing the hyperplanes whose independent sets correspond to the linearly independent subsets of the linear forms defining the hyperplanes.
This gives the matroids $M_0$ and $M$ associated to $\C_0$ and
$\widehat\C$, on $E_0$ and $E=E_0\sqcup\{0\}$, respectively.
The assumption that $\C$ is essential means that the forms $\ell_i$ span $(\R^r)^*$.
Thus $M_0$ has rank $r$, $M$ has rank $r+1$, and $M_0$ is the contraction of $M$ by $0$.

\subsection{Matroid conventions}
For a matroid $N$, write $\rk_N J$ for the rank of a subset $J$.
A \emph{flat} is a subset whose rank increases whenever an element
outside it is added. Let $N^J=N|J$ and $N_J=N/J$ denote restriction
and contraction, and let $N^*$ be the dual matroid. Write $\cL(N)$
for the lattice of flats and $\Z(N)$ for its cyclic flats, namely
the flats $J$ for which $N^J$ has no coloops. 

An affine matroid $(M,0)$ is a matroid on a ground set $E$ with a
choice of distinguished element $0\in E$. For such a pair, we write
\[
 \cL=\cL(M),\qquad \cL_0=\{F\in\cL\mid 0\in F\}.
\]
The flat correspondence for contraction identifies $\cL(M_0)$ with
$\cL_0$ by $K\mapsto K\cup\{0\}$. We distinguish these flats as
subsets of $E_0=E\setminus\{0\}$ and $E$, respectively.

For a loopless matroid $N$ of rank $s$, let
\[
 \chi_N(t)=\sum_{F\in\cL(N)}\mu(\varnothing,F)t^{s-\rk_N F},
\]
where $\mu$ is the M\"obius function of $\cL(N)$. We use
$\mup(N)=(-1)^s\chi_N(0)$ and, for $s>0$,
$\beta(N)=(-1)^{s-1}\chi_N'(1)$. For a matroid with loops we put
$\chi_N=0$ and $\mup(N)=\beta(N)=0$. For the empty matroid,
$\mup(\varnothing)=1$ and $\beta(\varnothing)=0$.

Returning to the affine matroid $(M,0)$, recall that for a flat $F \in \cL$ we set $w_F=\beta(M^F)\beta(M_F)$.

\subsection{Affine oriented matroids}
Central and affine arrangements are generalized by oriented matroids
and affine oriented matroids. These axiomatize the sign vectors of
faces, similarly to the way matroids axiomatize linear independence
among vectors.
An oriented matroid coming from a real hyperplane arrangement is
called \emph{realizable}. Not all oriented matroids are realizable.
Conjecturally, our determinant theorems extend to affine oriented matroids.
We recall the terminology.

An \emph{oriented matroid} $\M$ on a finite ground set $E$ is specified
by sign vectors $G\in\{0,+,-\}^E$, called \emph{covectors}, satisfying
the covector axioms \cite{OrientedMatroids1999}.
We write $M$ for its underlying matroid and assume that $M$ is simple
of rank $r+1$. An \emph{affine oriented matroid} $(\M,0)$ has a
distinguished element $0\in E$. A covector $G$ is \emph{feasible} if
$G(0)=+$. A \emph{tope} is a covector with no zero entries, and we
denote the feasible topes by $\cF$.

The homogenized arrangement $\widehat\C$ yields an affine oriented matroid $(\mathcal{M},0)$ on the ground set $E=E_0\cup\{0\}$ whose covectors are the sign vectors
\[
 \bigl(\operatorname{sign}\widehat f_i(x,t)\bigr)_{i\in E},
 \qquad (x,t)\in\R^{r+1}.
\]
The feasible covectors (topes) of $(\mathcal{M},0)$ correspond exactly to the sign vectors of all faces (chambers) in the affine arrangement $\C$.
Similarly, we can consider the oriented matroid $\mathcal M_0=\mathcal{M}/0$ which corresponds to the central arrangement $\C_0$.

Many combinatorial aspects of arrangements translate to oriented matroids.
The face order and zero set are given by
\[
 G\leq F\ \Longleftrightarrow\ G(i)\in\{0,F(i)\}\text{ for every }i,
 \qquad z(G)=\{i\mid G(i)=0\}.
\]
For every nonzero covector $G$, define
\[
 \dim G=r-\rk_M z(G).
\]
For feasible $G$, this agrees with the dimension of its affine face
in the realizable case. 

A feasible face $G$ is \emph{bounded} if
every nonzero covector $F\leq G$ is feasible. We denote the bounded
feasible topes by $\cB\subseteq\cF$. These correspond to the bounded
chambers of the affine arrangement.  As recalled in~\cite[Proposition 3.5]{Lam}, in both the case of hyperplane arrangements and more generally for affine oriented matroids, the number of bounded-regions (or bounded feasible topes) is given by
\[
 |\cB|=\beta(M).
\]

For an arrangement, write $H_J=\bigcap_{i\in J}H_i$.
We use this geometric notation for oriented matroids as well. 
 For $J\subseteq E_0$, write
$P^J$ for the sign vector obtained from $P$ by reversing the signs in
$J$, and abbreviate $P^{\{i\}}$ to $P^i$.

We also use the composition of sign vectors
\[
 (F\circ G)(i)=
 \begin{cases}F(i),&F(i)\ne0,\\G(i),&F(i)=0.\end{cases}
\]
Covectors are closed under composition. For an arrangement, this fills
the zero entries of $F$ according to $G$, describing a small displacement
from the face $F$ towards the face $G$.

We also use the \emph{strong elimination axiom} in the following
special case. If two covectors agree outside a coordinate $i$ and
have opposite signs at $i$, then setting that coordinate to zero
and keeping their common entries gives a covector.
For hyperplanes, this follows by intersecting the segment between
points of the two faces with $H_i$.

These definitions give the following dictionary.

\begin{center}
\small
\begin{tabularx}{\textwidth}{@{}XX@{}}
\toprule
Affine hyperplane arrangement & Affine oriented matroid\\
\midrule
Chamber $P\in\cF$ & Tope $P$ with $P(0)=+$\\
Face in $\bar P$ & Covector $G\leq P$ with $G(0)=+$\\
Face $G$ lies in the closure of $F$ & $G\leq F$ in the covector order\\
Hyperplanes containing $G$ & Zero set $z(G)$\\
Codimension of $G$ & $\rk_M z(G)$\\
Nonempty common face of $\bar P,\bar Q$ & Existence of a feasible $G\leq P,Q$\\
Hyperplanes separating $P,Q$ & $\Sep(P,Q)=\{i\mid P(i)\ne Q(i)\}$\\
$\bar P\cap H_J\ne\varnothing$ & Feasible $G\leq P$ with $J\subseteq z(G)$\\
Bounded face $G$ & Every nonzero $F\leq G$ has $F(0)=+$\\
Bounded chamber $P\in\cB$ & Bounded feasible tope\\
Small displacement from $F$ towards $G$ & Composition $F\circ G$\\
\bottomrule
\end{tabularx}
\end{center}

\section{Twisted intersection forms and period matrices}\label{sec:periods}
Our goal in this section and the next is to compute the determinant of 
$(V_b^{-1})^{\cB}$ in the realizable case.  In Section~\ref{subsec:betti-pairing} we recall Lam's result expressing the inverse of the Varchenko matrix in terms of his combinatorial Betti homology matrix.  Section~\ref{subsec:deRham-pairing} recalls Lam's combinatorial de Rham cohomology pairing and determinant formula for the pairing with respect to the weighted $\beta$nbc basis.  Section~\ref{subsec:geometric-comparison} summarizes Lam's identification of the combinatorial pairings with geometric pairings.  Section~\ref{subsec:period} introduces period matrices and Lam's version of the twisted period relation of Cho--Matsumoto.  In Section~\ref{subsec:period-determinant} we recall Varchenko, Loeser--Sabbah, Douai--Terao's determinant formula for the period matrices.  All of these ingredients are finally combined in Section~\ref{sec:determinants} to complete our proof of Theorem~\ref{thm:Lam-real}.

\subsection{Lam's combinatorial Betti pairing}\label{subsec:betti-pairing}
Let $(\M,0)$ be an affine oriented matroid of rank $r+1$ with ground set $E$, and let
$b_i$, $i\in E$, be independent variables, including $b_0$.
Lam's \emph{Betti homology matrix}
\[
 B_b=\bigl(B_b( P,Q) \bigr)_{P,Q\in\cF}
\]
is a symmetric matrix over $\QM(b_i\mid i\in E)$, defined by a sum over
flags of faces \cite[Definition 8.2]{Lam}. Its summands are signed
products of monomials $b_F$ and factors $(b_F^2-1)^{-1}$, where $F$
ranges over zero sets of nonzero faces of $P$. We will only need this
feature of the definition and the following formula for the inverse.

\begin{proposition}[{\cite[Theorem 9.2 and Corollary 9.3]{Lam}}]\label{prop:Betti-Varchenko}
For $P,Q\in\cF$,
\[
 (B_b^{-1})(P,Q)=\frac{(-1)^r}{1-b_E}
 \bigl(b_{\Sep(P,Q)}+(-1)^{r+1}b_{E\setminus\Sep(P,Q)}\bigr).
\]
In particular,
\[
 \left.B_b\right|_{b_0=0}=(-1)^rV_b^{-1}.
\]
\end{proposition}

We are interested in the $\cB \times \cB$-submatrix $B_b^{\cB}$.
This block is independent of $b_0$. Indeed, every nonzero face of a
bounded tope is feasible, so none of the zero sets in its face flags
contains $0$. Proposition~\ref{prop:Betti-Varchenko} therefore gives
\[
 B_b^{\cB}=\left( \left.B_b\right|_{b_0=0}\right)^\cB = (-1)^r(V_b^{-1})^{\cB}.
\]
Since $|\cB|=\beta(M)$, the determinants differ by the sign
\[
 \det B_b^{\cB}=(-1)^{r\beta(M)}
                 \det\bigl((V_b^{-1})^{\cB}\bigr).
\]
\begin{remark}
Theorem~\ref{thm:Lam-real} is stated in this Varchenko normalization.
The Betti determinant in \cite[Conjecture 8.10]{Lam} is missing this global sign, which is already visible
in \cite[Example 8.11]{Lam}.
\end{remark}
\begin{example}\label{ex:line}
For the two-point arrangement defined by $f_1(x)=x$ and $f_2(x)=1-x$,
the unique bounded chamber is $P=(0,1)$, oriented from $0$ to $1$.
All matrices in this example have size $1\times1$.
The Betti self-intersection of $P$ is
\[
 \det B_b^{\cB}=B_b(P,P)=1+\frac{1}{b_1^2-1}+\frac{1}{b_2^2-1}
     =-(V_b^{-1})(P,P) = -\det\bigl((V_b^{-1})^{\cB}\bigr),
\]
as in \cite[Example 14.2]{Lam}. Here the sign in
Proposition~\ref{prop:Betti-Varchenko} is $(-1)^r=-1$.
Example~\ref{ex:line4} recovers this entry from the
period relation.
\end{example}
\subsection{Lam's combinatorial de Rham pairing}\label{subsec:deRham-pairing}
Let $(M,0)$ be an affine matroid of rank $r+1$ with ground set $E = E_0 \sqcup \{0\}$, and let $a_i$, $i\in E$ be independent variables.
Let $A^{\bullet}(M)$ be the Orlik--Solomon algebra over $\CM$.
It is the exterior algebra on generators $e_i$, $i\in E$, modulo
the elements $\partial e_C$ for circuits $C$, where $\partial$ is
the derivation of degree $-1$ satisfying $\partial e_i=1$, and has top
degree $r+1$.
The reduced algebra is $\bar A^{\bullet}(M)=\ker\partial$.
It is generated by $\bar e_i=e_i-e_0$, $i\in E_0$, has top
degree $r$ and $A^{r+1}(M) \cong \bar A^r(M)$.  

Lam's combinatorial de Rham cohomology pairing is a bilinear form on
$\bar A^r(M)$, the top part of the reduced Orlik--Solomon algebra,
with values in rational functions in the variables $a_i$, $i\in E$:
\[
 \langle~,~\rangle^{\dR}:
 \bar A^r(M)\otimes_{\CM}\bar A^r(M)
 \longrightarrow\CM(a_i\mid i\in E).
\]
The bilinear form is defined by residues. For a flag of flats
\[
 F_\bullet=(\varnothing=F_0\subsetneq F_1\subsetneq\cdots
             \subsetneq F_r\subsetneq E),\qquad \rk F_j=j,
\]
let $\Res_{F_\bullet}:\bar A^r(M)\to\CM$ be the iterated
Orlik--Solomon residue of \cite[Definition 5.2]{Lam}. It takes residues
successively for the rank-one flats $F_j/F_{j-1}$ in $M/F_{j-1}$.
For $x,y\in\bar A^r(M)$, Lam's bilinear form is
\begin{equation}\label{eq:combinatorial-deRham}
 \langle x,y\rangle^{\dR}
 =\sum_{F_\bullet}
   \frac{\Res_{F_\bullet}(x)\Res_{F_\bullet}(y)}
        {a_{F_1}\cdots a_{F_r}},
\end{equation}
where the sum runs over all complete flags and $a_F = \sum_{i\in F} a_i$. This is \cite[Definition 5.6 and Corollary 5.12]{Lam} on the reduced
algebra.

This bilinear form descends to give a bilinear form on Aomoto cohomology.  To define Aomoto cohomology, we specialize the variables $a_i, i \in E$ to complex numbers such that $\sum_{i\in E} a_i = 0$.  The Aomoto cohomology $A_{\mathbf a}$ depends on the choice $\mathbf a=(a_i)_{i\in E}$ of specialization.  For any $J \subseteq E$, let $a_J=\sum_{i\in J}a_i$.  In addition to the assumption that $a_E=\sum_{i\in E} a_i =0$ (equivalently, $a_0=-\sum_{i\in E_0}a_i$) we assume that the weights $a_i, i\in E$ are \emph{nonresonant}, meaning
\begin{equation}\label{eq:weights}
 a_F\notin\ZM\quad
 \text{ for all }F\in\cL\setminus\{\varnothing,E\}\text{ with }\ M^F\text{ connected}.
\end{equation}
While the denominators in the expression for $\langle x,y\rangle^{\dR}$ above are products ${a_{F_1}\cdots a_{F_r}}$ for arbitrary flags, Lam~\cite[Corollary 10.25]{Lam} shows that $\langle x,y\rangle^{\dR}$ can be re-expressed such that only products of $a_F$ for \emph{connected} flats $F \neq \hat 0, \hat 1$ appear in the denominator.  In particular, under the condition (\ref{eq:weights}), the rational function $\langle x,y\rangle^{\dR}$ on $a_E=0$ is regular.
The nonresonance condition is also used in the basis statements below and the geometric comparisons in
Section~\ref{subsec:geometric-comparison}.

Put $\omega_{\mathbf a}=\sum_{i\in E_0}a_i\bar e_i$.
Equipping $\bar A^{\bullet}(M)$ with differential
$\omega_{\mathbf a}\wedge(-)$ gives the \emph{Aomoto complex}.
The \emph{Aomoto cohomology} is defined as
\[
 A_{\mathbf a}
 :=H^r(\bar A^{\bullet}(M),\omega_{\mathbf a}\wedge(-))
 =\frac{\bar A^r(M)}{\omega_{\mathbf a}\wedge\bar A^{r-1}(M)}.
\]
Thus $A_{\mathbf a}=A_{-\mathbf a}$ as quotients of $\bar A^r(M)$.
We retain the two labels to distinguish their later geometric
realizations with dual local systems.
Under \eqref{eq:weights}, the Aomoto cohomology is concentrated in
degree $r$ and has dimension $\beta(M)$ \cite[Theorem 7.3]{Lam}.

At nonresonant weights, Lam's de Rham cohomology pairing $\langle~,~\rangle^{\dR}$ descends to a nondegenerate symmetric pairing on $A_{\mathbf a}$, which we denote by $\langle~,~\rangle^{\dR}_{\mathbf a}$ 
\cite[Theorem 7.7 and Corollary 7.15]{Lam}.

We now recall the \textit{weighted $\beta$nbc basis} of $A_{\mathbf a}$. Fix a total order on $E$ with $0$ least.
A broken circuit is a circuit
with its least element removed. An nbc basis contains no broken
circuit. It is a $\beta$nbc basis if for every $i\in B\setminus\{0\}$
there is $j<i$ such that $(B\setminus\{i\})\cup\{j\}$ is again a
basis. Write $\Bnbc$ for this set. Its cardinality is $\beta(M)$
\cite{Ziegler1992}, \cite[Definition 7.9 and Theorem 7.10]{Lam}.

Write a $\beta$nbc basis in decreasing order as
$B=(i_1>\cdots>i_r>i_{r+1}=0)$ and put
\[
 F_j(B)=\cl_M\{i_1,\ldots,i_j\},\qquad
 \omega_{\mathbf a}(F)=\sum_{i\in F}a_i\bar e_i.
\]
Here $0\notin F_j(B)$. The corresponding weighted form is
\begin{equation}\label{eq:nbc-form}
 \psi_{\mathbf a}^B
 =\omega_{\mathbf a}(F_r(B))\wedge\cdots\wedge
   \omega_{\mathbf a}(F_1(B)).
\end{equation}
At nonresonant weights these forms give a basis of $A_{\mathbf a}$
\cite{FalkTerao1997}, \cite[Corollary 7.14]{Lam}.
They satisfy $\psi_{-\mathbf a}^B=(-1)^r\psi_{\mathbf a}^B$.

Set
\[
 D_{\mathbf a}^{\mathrm{nbc}}(B,B')
   =\langle\psi_{\mathbf a}^B,\psi_{-\mathbf a}^{B'}\rangle^{\dR}_{\mathbf a}.
\]
This uses the basis at its own weight in each argument. Relative
to Lam's Gram matrix using $\psi_{\mathbf a}$ in both arguments,
$D_{\mathbf a}^{\mathrm{nbc}}$ is multiplied by $(-1)^r$.

Recall $w_F=\beta(M^F)\beta(M_F)$ and put
\[
 R_M(\mathbf a)
 =\frac{\displaystyle\prod_{F\in\cL_0\setminus\{E\}}a_F^{w_F}}
        {\displaystyle\prod_{F\in\cL\setminus(\cL_0\cup\{\varnothing\})}
             a_F^{w_F}}.
\]

\begin{theorem}[Lam]\label{thm:deRham-determinant}
As rational functions on $a_E=0$,
\[
 \det D_{\mathbf a}^{\mathrm{nbc}}=\pm R_M(\mathbf a)^{-1}.
\]
The sign is constant.
\end{theorem}
\begin{proof}
This is \cite[Corollary 7.20]{Lam}.
The use of $\psi_{-\mathbf a}$ in the second argument changes the
determinant by the constant sign $(-1)^{r\beta(M)}$.
\end{proof}

\begin{example}\label{ex:line2}
Consider again the two-point arrangement from Example~\ref{ex:line}.
Its matroid $M$ is the uniform matroid $U_{2,3}$ on $E=\{0,1,2\}$.
For the order $0<1<2$, the nbc bases are $\{0,1\}$ and $\{0,2\}$,
and the unique $\beta$nbc basis is $B=\{0,2\}$.
The weighted form is $\psi^B_{\mathbf a}=a_2\bar e_2$.
Since $\psi^B_{-\mathbf a}=-\psi^B_{\mathbf a}$,
\eqref{eq:combinatorial-deRham} and $a_0=-a_1-a_2$ give
\[
 \det D_{\mathbf a}^{\mathrm{nbc}} = D_{\mathbf a}^{\mathrm{nbc}}(B,B)
 =-a_2^2\left(\frac1{a_2}+\frac1{a_0}\right)
 =\frac{a_1a_2}{a_0}.
\]
Note that this agrees with Lam's Theorem~\ref{thm:deRham-determinant}.
\end{example}

\begin{remark}\label{rem:chamber-deRham}
For an orientation $\M$ of $M$, Eur and Lam associate to each tope $P$ a weight-independent canonical
element $\Omega_P\in A^{r+1}(M)$ \cite[Theorem 2.11]{EurLam}.
Put $\bar\Omega_P=\partial\Omega_P$.
At nonresonant weights their classes for $P\in\cB$ form a second
basis of $A_{\mathbf a}$ \cite[Corollary 7.14]{Lam}. Its Gram matrix
\[
 D_{\mathbf a}^{\mathrm{ch}}(P,Q)
 =\langle\bar\Omega_P,\bar\Omega_Q\rangle^{\dR}_{\mathbf a}
\]
has determinant $\pm R_M(\mathbf a)$ by
\cite[Theorem 7.18]{Lam}.

However, we focus our attention on the Gram matrix for the weighted $\beta$nbc basis (instead of the canonical chamber basis) because, as recalled in Section~\ref{subsec:period-determinant}, the determinants of the corresponding period matrices are known.
\end{remark}

\subsection{Geometric realizations and comparison}\label{subsec:geometric-comparison}
We will now recall Lam's geometric realizations of the Betti and de Rham matrices defined above. For this, in this section we assume that $\C$ is a real essential affine arrangement
and the weights $a_i$ are nonresonant, that is, \eqref{eq:weights} holds.
Let
\[
 U=\CM^r\setminus\bigcup_{i\in E_0}H_i^{\CM}
\]
be the complement of the complexified arrangement. Specialize
$b_i=e^{-\pi i a_i}$, so $b_E=1$. We use the bounded block
$B_b^{\cB}$, which is independent of $b_0$.

The logarithmic connection
\[
 \nabla_{\mathbf a}=\dd+\omega_{\mathbf a}\wedge(-),\qquad
 \omega_{\mathbf a}=\sum_{i\in E_0}a_i\,\dd\log f_i
\]
has a rank-one local system $\Loc_{\mathbf a}$ of flat sections.
Locally these are scalar multiples of
$\varphi_{\mathbf a}^{-1}=\prod_{i\in E_0}f_i^{-a_i}$.
The monodromy around $H_i$ is $b_i^2$ and
$\Loc_{\mathbf a}^{\vee}=\Loc_{-\mathbf a}$.

We use ordinary and Borel--Moore homology, and ordinary and compactly
supported cohomology, with coefficients in these local systems.
Under \eqref{eq:weights}, the natural maps
\[
 H_r(U,\Loc_{\pm\mathbf a})\longrightarrow
 H_r^{\BM}(U,\Loc_{\pm\mathbf a}),\qquad
 H_c^r(U,\Loc_{\pm\mathbf a})\longrightarrow
 H^r(U,\Loc_{\pm\mathbf a})
\]
are isomorphisms by \cite[Theorem 12.3]{Lam}. Kohno proved this under
stronger nonresonance assumptions in \cite[Theorem 1]{Kohno}.
We denote their inverses by $\reg$.
The natural maps, and hence their inverses, are dual under the perfect
evaluation pairings \cite[equation (12.3)]{Lam}. Explicitly,
\[
 (\xi,\reg c)=(\reg\xi,c),
 \qquad \xi\in H^r(U,\Loc_{-\mathbf a}),\quad
 c\in H_r^{\BM}(U,\Loc_{\mathbf a}).
\]
The left side pairs cohomology with ordinary homology, and the right
side pairs compactly supported cohomology with Borel--Moore homology.

These groups vanish outside degree $r$, where they have dimension
$\beta(M)$ \cite[Theorem 12.3 and Section 12.3]{Lam}.
The algebraic de Rham groups are
\[
 H^\bullet_{\dR}(U,\nabla_{\pm\mathbf a})
 =H^\bullet(\Gamma(U,\Omega_U^\bullet),\nabla_{\pm\mathbf a}).
\]
The de Rham comparison gives isomorphisms
\[
 \comp_{\pm\mathbf a}:H^\bullet_{\dR}(U,\nabla_{\pm\mathbf a})
 \xrightarrow{\ \sim\ }H^\bullet(U,\Loc_{\pm\mathbf a})
\]
\cite[Section 12.3]{Lam}. We also write 
\[\reg: H^r_{\dR}(U,\nabla_{\pm\mathbf a}) \to H^r_{c,\dR}(U,\nabla_{\pm\mathbf a}) \] for the corresponding
regularization in de Rham cohomology.

Consider the logarithmic realization
\[
 \iota:\bar A^{\bullet}(M)\hookrightarrow\Gamma(U,\Omega_U^{\bullet}),
 \qquad \bar e_i\longmapsto\dd\log f_i=\frac{\dd f_i}{f_i},
\]
which induces a map of complexes
\[
 (\bar A^{\bullet}(M),\omega_{\mathbf a}\wedge(-))
 \xhookrightarrow{\ \iota\ }
 (\Gamma(U,\Omega_U^{\bullet}),\nabla_{\mathbf a}).
\]
Under \eqref{eq:weights}, this is a quasi-isomorphism by the theorem
of Esnault--Schechtman--Viehweg \cite{EsnaultSchechtmanViehweg1992}
and Schechtman--Terao--Varchenko \cite{SchechtmanTeraoVarchenko1995},
as stated in \cite[Theorem 12.4]{Lam}.

Taking top cohomology gives the identifications
\[
 A_{\pm\mathbf a}
 \xrightarrow{\ \iota_*\ }H^r_{\dR}(U,\nabla_{\pm\mathbf a})
 \xrightarrow{\ \comp_{\pm\mathbf a}\ }H^r(U,\Loc_{\pm\mathbf a}).
\]
The first sends $[\eta]$ to the twisted de Rham class $[\iota(\eta)]$.
The second is the de Rham comparison with the local system of flat
sections.
Write $\per_{\pm\mathbf a}=\comp_{\pm\mathbf a}\circ\iota_*$
for the composite.
Under $\iota$, the elements $\bar\Omega_P$ become the canonical
rational top forms of the chamber closures, with logarithmic poles
along their boundaries \cite[Remark 4.10]{Lam}.

Choose once and for all the branch $\log$ with cut along the positive
imaginary ray, agreeing with the real logarithm on $\R_{>0}$. We write
\[
 z^a=\exp(a\log z).
\]
The product
\[
 \varphi_{\mathbf a}
 =\prod_{i\in E_0}f_i^{a_i}
\]
is holomorphic on the complex neighbourhood $P+i\R^r$ of every real
chamber $P$.
Put $P^-=\{i\in E_0\mid P(i)=-\}$ and define the normalized branches by
\begin{equation}\label{eq:loading}
 \varphi_{\mathbf a,P}
 =b_{P^-}^{-1}\varphi_{\mathbf a},\qquad
 \varphi_{\mathbf a,P}^{-1}=b_{P^-}\varphi_{\mathbf a}^{-1}.
\end{equation}
Its restriction to the real chamber satisfies
\[
 \varphi_{\mathbf a,P}(x)=\prod_{i\in E_0}|f_i(x)|^{a_i},
 \qquad x\in P.
\]
Thus it is positive on $P$ for real weights $\mathbf a$, and
$\varphi_{-\mathbf a,P}=\varphi_{\mathbf a,P}^{-1}$.
The chosen branch extends by analytic continuation along paths in $U$.
This is the normalization of loaded chambers in
\cite[Definition 12.6 and equation (12.5)]{Lam}.

With the ambient orientation of $\R^r$, the loaded chambers define
classes in Borel--Moore homology
\begin{align*}
 C_{\mathbf a}(P)&=[P\otimes\varphi_{\mathbf a,P}^{-1}]
       \in H_r^{\BM}(U,\Loc_{\mathbf a}),\\
 C_{-\mathbf a}(P)&=[P\otimes\varphi_{\mathbf a,P}]
       \in H_r^{\BM}(U,\Loc_{-\mathbf a}).
\end{align*}
The classes $C_{\mathbf a}(P)$ and $C_{-\mathbf a}(P)$, for $P\in\cB$, are bases
of the respective Borel--Moore groups
\cite{Kohno}, \cite[Proposition 12.7]{Lam}.

The Betti homology pairing
\[
 \langle-,-\rangle_{\mathrm{hom}}:
 H_r^{\BM}(U,\Loc_{-\mathbf a})\otimes
 H_r^{\BM}(U,\Loc_{\mathbf a})\longrightarrow\CM
\]
is obtained by applying $\reg$ to the second argument and then the
Poincar\'e--Verdier duality pairing \cite[Definition 12.5]{Lam}.
Using regularization to identify Borel--Moore and ordinary homology,
its dual under the perfect evaluation pairings is the Betti
cohomology pairing
\[
 \langle-,-\rangle_{\mathrm{coh}}:
 H^r(U,\Loc_{\mathbf a})\otimes H^r(U,\Loc_{-\mathbf a})
 \longrightarrow\CM
\]
\cite[equation (12.3) and Section 12.4]{Lam}.

Similarly, we obtain a pairing on de Rham cohomology which is normalized by $(2\pi i)^{-r}$. For a class
$x$ at weight $\mathbf a$ and a class $y$ at weight $-\mathbf a$,
we use the argument order of \cite[equation (15.1)]{Lam} and write
\[
 \langle x,y\rangle^{\nabla}_{\mathbf a}
 =\frac{1}{(2\pi i)^r}\int_U y\wedge\reg(x).
\]
Here the integral uses a compactly supported smooth representative
of $\reg(x)$.

Lam's comparison theorems~\cite[Theorems 12.9 and 12.12]{Lam} identify the geometric pairings with Lam's combinatorial
Betti and de Rham pairings:
\begin{align*}
 \langle C_{-\mathbf a}(P),C_{\mathbf a}(Q)\rangle_{\mathrm{hom}}
   &=B_b^{\cB}(P,Q),\\
 \langle\iota_*x,\iota_*y\rangle^{\nabla}_{\mathbf a}
   &=\langle x,y\rangle^{\dR}_{\mathbf a}.
\end{align*}
Here $P,Q\in\cB$, $x\in A_{\mathbf a}$ and $y\in A_{-\mathbf a}$.
 In particular,
$D_{\mathbf a}^{\mathrm{nbc}}$ and $D_{\mathbf a}^{\mathrm{ch}}$ are also the Gram
matrices of the normalized geometric de Rham pairing in the two
displayed bases.

Under the perfect evaluation pairings $(~,~)$, let
$\{\xi_{\pm\mathbf a}(P)\}_{P\in\cB}$ be the bases of
$H^r(U,\Loc_{\pm\mathbf a})$ dual to
$\{\reg C_{\mp\mathbf a}(P)\}_{P\in\cB}$, respectively.
Passing to dual bases takes the inverse transpose Gram
matrix \cite[Section 12.4]{Lam}, so symmetry of $B_b^{\cB}$ gives
\begin{equation}\label{eq:dual-Gram}
 \langle\xi_{\mathbf a}(P),\xi_{-\mathbf a}(Q)\rangle_{\mathrm{coh}}
 =\bigl((B_b^{\cB})^{-1}\bigr)(P,Q).
\end{equation}

\subsection{The period relation}\label{subsec:period}
For a real arrangement and nonresonant weights, define the regularized
$\cB\times\Bnbc$ period matrix by
\[
 \PM^*_{\mathbf a}(P,B)
 =\bigl(\per_{\mathbf a}(\psi_{\mathbf a}^B),\reg C_{-\mathbf a}(P)\bigr).
\]
Thus rows of the period matrix are indexed by bounded chambers and columns by $\beta$nbc bases of $M$.
The parentheses denote the evaluation pairing between cohomology
and homology with dual local coefficients. This is the notation of
\cite[Definition 3.3.2]{DouaiTerao1997}.
Since the $\xi_{\mathbf a}(P)$ are evaluation-dual to the regularized
chamber cycles, we have
\begin{equation}\label{eq:period-coordinates}
 \per_{\mathbf a}(\psi_{\mathbf a}^B)
 =\sum_{P\in\cB}\PM^*_{\mathbf a}(P,B)\,\xi_{\mathbf a}(P).
\end{equation}

For $\operatorname{Re}a_i>0$, $i\in E_0$, the chamber integrals
converge and define the matrix
\[
 \PM_{\mathbf a}(P,B)
 =\int_P\varphi_{\mathbf a,P}\,\iota(\psi_{\mathbf a}^B).
\]
At nonresonant weights in this region, $\PM_{\mathbf a}=\PM^*_{\mathbf a}$
\cite[Remark 3.3.3]{DouaiTerao1997}.
The matrix $\PM^*_{-\mathbf a}$ uses the opposite loading
and the forms $\psi_{-\mathbf a}^B$.

\begin{example}\label{ex:line3}
Continuing Examples~\ref{ex:line} and~\ref{ex:line2}, the form
$\psi^B_{\mathbf a}=a_2\bar e_2$ is realized as
$\iota(\psi^B_{\mathbf a})=-a_2\,\dd x/(1-x)$.
For $\operatorname{Re}a_1,\operatorname{Re}a_2>0$, its period is
\[
 \PM_{\mathbf a}(P,B)
 =-a_2\int_0^1 x^{a_1}(1-x)^{a_2-1}\,\dd x
 =-a_2 B(1+a_1,a_2),
\]
where $B(x,y)$ is the classical beta function. The expression
\[
 \PM_{\mathbf a}^*(P,B)
 =-\frac{\Gamma(1+a_1)\Gamma(1+a_2)}{\Gamma(1-a_0)}
\]
gives its meromorphic continuation.
\end{example}

\begin{proposition}\label{prop:period-relation}
For nonresonant weights with $b_i=e^{-\pi i a_i}$,
\begin{equation}\label{eq:period-relation}
 (2\pi i)^rD_{\mathbf a}^{\mathrm{nbc}}
 =(\PM^*_{\mathbf a})^{\mathsf T}(B_b^{\cB})^{-1}\PM^*_{-\mathbf a}.
\end{equation}
Consequently,
\begin{equation}\label{eq:det-period-relation}
 \det B_b^{\cB}=(2\pi i)^{-r\beta(M)}
       \frac{\det\PM^*_{\mathbf a}\det\PM^*_{-\mathbf a}}
            {\det D_{\mathbf a}^{\mathrm{nbc}}}.
\end{equation}
\end{proposition}
\begin{proof}
The comparison theorem gives
\[
 (2\pi i)^r\langle x,y\rangle^{\dR}_{\mathbf a}
 =\langle\per_{\mathbf a}(x),\per_{-\mathbf a}(y)\rangle_{\mathrm{coh}}
\]
by \cite[equation (15.1)]{Lam}. Substituting
\eqref{eq:period-coordinates} at both weights and using
\eqref{eq:dual-Gram} gives \eqref{eq:period-relation}.
This is the twisted period relation of Cho--Matsumoto
\cite{ChoMatsumoto1995} in the form of
\cite[equation (15.2)]{Lam}. Taking determinants gives
\eqref{eq:det-period-relation}.
\end{proof}

\subsection{Period determinants and continuation}\label{subsec:period-determinant}
For a bounded region $P$ and $i\in E_0$, let $\Sigma_{i,P}$ be
the face of $\bar P$ on which $|f_i|$ is maximal. Set
\[
 c_{i,P}=|f_i(\Sigma_{i,P})|>0,\qquad
 C(\mathbf a)=\prod_{i\in E_0}\prod_{P\in\cB}c_{i,P}^{a_i}.
\]
These are the external supporting-face critical values in
\cite[Sections 2.6 and 4.3.2]{MukherjeeVarchenko2009}, expressed in
our positive normalization. In particular,
\begin{equation}\label{eq:critical-cancellation}
 C(\mathbf a)C(-\mathbf a)=1.
\end{equation}
The beta function of the weighted arrangement is
\[
 \Betafunc_M(\mathbf a)
 =\frac{\displaystyle
      \prod_{F\in\cL\setminus(\cL_0\cup\{\varnothing\})}
          \Gamma(1+a_F)^{w_F}}
        {\displaystyle
      \prod_{F\in\cL_0\setminus\{E\}}\Gamma(1-a_F)^{w_F}}.
\]

\begin{theorem}[Loeser--Sabbah, Douai--Terao]\label{thm:period-determinant}
For $\operatorname{Re}a_i>0$, $i\in E_0$,
\[
 \det\PM_{\mathbf a}=\varepsilon\, C(\mathbf a)\Betafunc_M(\mathbf a),
 \qquad\varepsilon\in\{1,-1\},
\]
where $\varepsilon$ depends only on the fixed orderings and orientations.
\end{theorem}
\begin{proof}
We use \cite{LoeserSabbah1991} and
\cite[Theorem 5.1.1]{DouaiTerao1997} in the formulation of
\cite[Theorem 4.1 and Sections 2.5, 4.2.1, 4.3]{MukherjeeVarchenko2009}.
The forms \eqref{eq:nbc-form} are their weighted flag forms, the
exponents are $w_F$, and the critical-value product for the loading
\eqref{eq:loading} is $C(\mathbf a)$. Changing the orderings or
orientations contributes the fixed sign $\varepsilon$.
\end{proof}

The entries of $\PM_{\mathbf a}$ are holomorphic for
$\operatorname{Re}a_i>0$, $i\in E_0$. Their meromorphic continuations
to the full weight space $a_E=0$ are the entries of $\PM^*_{\mathbf a}$
\cite[Remark 3.3.3]{DouaiTerao1997}.
Since the determinant is a polynomial in the entries,
Theorem~\ref{thm:period-determinant} therefore gives
\[
 \det\PM^*_{\mathbf a}=\varepsilon\,C(\mathbf a)\Betafunc_M(\mathbf a)
\]
as an identity of meromorphic functions on $a_E=0$. We can now
apply it at both $\mathbf a$ and $-\mathbf a$, using the geometric
comparisons at nonresonant weights.

\begin{corollary}\label{cor:opposite-periods}
As meromorphic functions,
\[
 \det\PM^*_{\mathbf a}\det\PM^*_{-\mathbf a}
 =\Betafunc_M(\mathbf a)\Betafunc_M(-\mathbf a).
\]
\end{corollary}
\begin{proof}
Multiply the continued formulas at $\mathbf a$ and $-\mathbf a$.
The critical-value factors cancel by \eqref{eq:critical-cancellation},
and the signs multiply to $\varepsilon^2=1$.
\end{proof}

\begin{remark}
Using $\varphi_{\mathbf a}$ in place of the normalized branches in
\eqref{eq:loading} replaces $B_b^{\cB}$ by
$\Lambda B_b^{\cB}\Lambda^{-1}$, where
$\Lambda=\diag(b_{P^-})_{P\in\cB}$.
The period matrices $\PM^*_{\mathbf a}$ and $\PM^*_{-\mathbf a}$
are multiplied on the left by $\Lambda$ and $\Lambda^{-1}$,
respectively. Thus the Betti determinant and the product of period
determinants are unchanged.
\end{remark}

\section{Proof of Lam's determinant conjecture}\label{sec:determinants}

\begin{proof}[Proof of Theorem~\ref{thm:Lam-real}]

Theorem~\ref{thm:deRham-determinant},
Corollary~\ref{cor:opposite-periods}, and
\eqref{eq:det-period-relation} give
\[
 \det B_b^{\cB}
 =\pm(2\pi i)^{-r\beta(M)}
       R_M(\mathbf a)\Betafunc_M(\mathbf a)\Betafunc_M(-\mathbf a).
\]
Substituting the products defining $R_M$ and $\Betafunc_M$ gives
\[
 \det B_b^{\cB}
 =\pm(2\pi i)^{-r\beta(M)}
 \frac{\displaystyle
   \prod_{F\in\cL\setminus(\cL_0\cup\{\varnothing\})}
      \bigl(a_F^{-1}\Gamma(1+a_F)\Gamma(1-a_F)\bigr)^{w_F}}
      {\displaystyle
   \prod_{F\in\cL_0\setminus\{E\}}
      \bigl(a_F^{-1}\Gamma(1+a_F)\Gamma(1-a_F)\bigr)^{w_F}}.
\]
The reflection formula and $b_F=e^{-\pi i a_F}$ give
\begin{equation}\label{eq:reflection}
 a_F^{-1}\Gamma(1+a_F)\Gamma(1-a_F)
 =\frac{\pi}{\sin(\pi a_F)}
 =\frac{2\pi i b_F}{1-b_F^2}.
\end{equation}
Collecting the powers of $2\pi i$ gives
\[
 \det B_b^{\cB}
 =\pm(2\pi i)^\delta\frac{\displaystyle
      \prod_{F\in\cL\setminus(\cL_0\cup\{\varnothing\})}
          \bigl(b_F/(1-b_F^2)\bigr)^{w_F}}
        {\displaystyle
      \prod_{F\in\cL_0\setminus\{E\}}
          \bigl(b_F/(1-b_F^2)\bigr)^{w_F}}.
\]
Here $\delta$ is the total exponent of $2\pi i$.
For $F\in\cL_0$, the relation $b_E=1$ gives
\[
 1-b_F^2=-b_F^2(1-b_{E\setminus F}^2).
\]
Using Proposition~\ref{prop:Betti-Varchenko} and absorbing the signs
into $c=\pm(2\pi i)^\delta$ gives
\[
 \det\bigl((V_b^{-1})^{\cB}\bigr)
 =c\left(\prod_{\varnothing\ne F\ne E}b_F^{w_F}\right)
 \frac{\displaystyle\prod_{F\in\cL_0\setminus\{E\}}
            (1-b_{E\setminus F}^2)^{w_F}}
      {\displaystyle\prod_{F\in\cL\setminus(\cL_0\cup\{\varnothing\})}
            (1-b_F^2)^{w_F}}.
\]
By \cite[Lemma 2.15]{LopezdeMedranoRinconShaw2020}, the sum
\[
 s=\sum_{\substack{F\in\cL\setminus\{\varnothing,E\}\\i\in F}}w_F
\]
is independent of $i\in E$. This also follows from
\cite[Corollary 2.1.15]{Gioan2002}. Hence the monomial factor is
\[
 \prod_{\varnothing\ne F\ne E}b_F^{w_F}
 =\prod_{i\in E}b_i^s=b_E^s=1.
\]
This proves \eqref{eq:Lam-determinant} up to the constant $c$. After eliminating
$b_0$, both sides are rational functions of the finite parameters
$b_i$, $i\in E_0$. Their equality on nonresonant weights therefore
gives an identity of rational functions. Both sides are regular at
$b_i=0$ for all $i\in E_0$, where $V_b$ becomes the identity and the
product becomes $1$. Thus $c=1$, proving the theorem.
\end{proof}

We obtain the following identity as a corollary.
\begin{corollary}\label{cor:degree}
For an affine oriented matroid $(\M,0)$ of rank $r+1$,
\[
 \sum_{F\in\cL\setminus(\cL_0\cup\{\varnothing\})}w_F
 -\sum_{F\in\cL_0\setminus\{E\}}w_F=r\beta(M).
\]
\end{corollary}
\begin{proof}
For real arrangements, $\delta$ in the preceding proof is the
left-hand side minus $r\beta(M)$. Since $1=|c|=(2\pi)^\delta$,
we have $\delta=0$.
For arbitrary affine oriented matroids,
\eqref{eq:combinatorial-deRham} and \eqref{eq:nbc-form} show that
each entry of $D_{\mathbf a}^{\mathrm{nbc}}$ is homogeneous of degree
$r$ in the weights. Its determinant has degree $r\beta(M)$. Comparing with
$R_M(\mathbf a)^{-1}$ in Theorem~\ref{thm:deRham-determinant}
gives the same identity.
\end{proof}

\begin{example}\label{ex:line4}
We conclude by computing the period relation for the two-point
arrangement. Substituting Examples~\ref{ex:line2} and~\ref{ex:line3}
into \eqref{eq:det-period-relation} and applying~(\ref{eq:reflection}) gives
\[
 \begin{aligned}
 \det B_b^{\cB} = B_b(P,P)
 &=\frac{a_0}{2\pi i a_1a_2}
   \frac{\Gamma(1+a_1)\Gamma(1-a_1)
         \Gamma(1+a_2)\Gamma(1-a_2)}
        {\Gamma(1+a_0)\Gamma(1-a_0)}\\
 &=\frac{b_1 b_2(1-b_0^2)}{b_0 (1-b_1^2)(1-b_2^2)}.
 \end{aligned}
\]
Since $b_0b_1b_2=1$,
\[
 \begin{aligned}
 \det B_b^{\cB}
 &=\frac{(1-b_0^2)}{b_0^2 (1-b_1^2)(1-b_2^2)}\\
 &=-\frac{1-b_1^2b_2^2}{(1-b_1^2)(1-b_2^2)}
 =-(V_b^{-1})^{\cB},
 \end{aligned}
\]
and we recover the Betti entry in Example~\ref{ex:line}.
\end{example}

\section{Combinatorial lemmas}\label{sec:comb}

In the remainder of the paper we will show that Theorem~\ref{thm:Lam-real} implies Theorem~\ref{thm:intersection-real} and more generally for affine oriented matroids, equation~(\ref{eq:Lam-determinant}) implies equation~(\ref{eq:Sq-determinant}) in the case of normal crossings.  In this section we assemble the necessary combinatorial definitions and properties of normal crossing hyperplane arrangements and affine oriented matroids.

\subsection{Normal crossing} We will make use of the following important genericity condition.
\begin{definition}\label{def:normal-crossing}
An affine oriented matroid $(\M,0)$ is \emph{normal crossing} if the
zero set $z(G)$ of every feasible covector $G$ is independent in the
underlying matroid $M$.
\end{definition}
For an affine hyperplane arrangement, this means that every nonempty
intersection of $s$ hyperplanes has codimension $s$. Equivalently, the
reduced divisor
\[
 D=\bigcup_{i\in E_0}H_i^{\CM}\subseteq\CM^r
\]
is a simple normal crossing divisor.
This condition holds for generic parameters $c_i$ in the equations $f_i(x)=\ell_i(x)+c_i$.

In oriented-matroid terminology, it corresponds to the genericity of the distinguished
element $0$, see \cite[Remark 2.3 and Section 5.2]{KowalenkoMautner}.
Normal crossing makes the local face structure Boolean, which we use
in the quadratic presentations and the Koszul complex below.

Normal crossing also implies that $M$ is the unique free coextension of $M_0$, the contraction of $M$ by $0$. This explains why
the normal crossing determinant formula involves only $M_0$.

\begin{proposition}\label{prop:normal-crossing-ranks}
If $(\M,0)$ is normal crossing, then $M$ is the free coextension of $M_0$ by $0$.
For $A\subseteq E_0$, its rank function is
\begin{equation}\label{eq:coextension-rank}
 \rk_M(A\cup\{0\})=\rk_{M_0}A+1,\qquad
 \rk_M A=\min\{\rk_{M_0}A+1,|A|\}.
\end{equation}
Its proper flats avoiding $0$ are the independent subsets of $M_0$.
Its proper flats containing $0$ are $K\cup\{0\}$ for proper flats
$K$ of $M_0$.
\end{proposition}
\begin{proof}
Every flat avoiding $0$ is the zero set of a feasible covector,
after reversing all signs if necessary. Normal crossing makes such
a flat independent. If $A$ is independent in $M_0$, it is independent
in $M$. If $A$ is dependent in $M_0$, then $0\in\cl_M A$, since
otherwise $\cl_M A$ would be an independent flat avoiding $0$.
The contraction rank formula gives \eqref{eq:coextension-rank}.
This is the rank function of the free coextension
\cite[Sections 7.2--7.3]{Oxley2011}, and the flat description follows.
\end{proof}

To compute the weights $w_F$ when $(\M,0)$ is normal crossing, we
use the following identity for the beta invariant.
For a matroid $N$ on $E_0$, write $\widehat N$ for its free coextension
by $0$. Its rank function is given by \eqref{eq:coextension-rank}, with
$(M,M_0)$ replaced by $(\widehat N,N)$.

\begin{lemma}\label{lem:coextension-beta}
For any matroid $N$,
\[
 \beta(\widehat N)=\mup(N^*).
\]
\end{lemma}
\begin{proof}
This follows from Brylawski's free extension formula for the beta
invariant and duality
\cite[Remarks 6.15(3)(b) and (d), p.~241]{Brylawski1982}.
\end{proof}

\begin{lemma}\label{lem:free-coextension}
Let $(\M,0)$ be normal crossing and $M_0=M/0$. For a nonempty proper
flat $F$ of $M$, 
\[
 w_F= \beta(M^F)\beta(M_F) = 
 \begin{cases}
  \beta(M_{\{i\}}),&F=\{i\},\quad i\in E_0,\\
  0,&0\notin F,\quad |F|>1,\\
  \beta(M_{0,K})\mup((M_0^K)^*),&F=K\cup\{0\}.
 \end{cases}
\]
\end{lemma}
\begin{proof}
By Proposition~\ref{prop:normal-crossing-ranks}, for a flat avoiding
$0$, $M^F$ is free and its beta invariant vanishes unless
$F$ is a singleton. For $F=K\cup\{0\}$,
\[
 M^F=\widehat{M_0^K},\qquad M_F=M_{0,K}.
\]
Apply Lemma~\ref{lem:coextension-beta} to the restriction.
\end{proof}

\subsection{Face lattices} The following are standard observations for affine hyperplane
arrangements. We use the covector axioms and the cited results to
obtain their oriented-matroid versions.

We write
$\bar P\cap\bar Q\cap H_J\ne\varnothing$ if there is a
feasible $G\leq P,Q$ with $J\subseteq z(G)$. Moreover, we define
\[
 \Delta_{PQ}=\{J\subseteq E_0\mid\bar P\cap\bar Q\cap H_J\ne\varnothing\}.
\]
Thus $\Delta_{PQ}$ records the collections of hyperplanes whose
intersection meets the common face of $P$ and $Q$.
It is closed under taking subsets.

Near a normal crossing face, either side of each hyperplane containing
the face can be chosen independently. This gives a cube of feasible
regions. Conversely, if all sign choices on $S$ occur while the signs
outside $S$ are fixed, convexity gives a common face on the hyperplanes
in $S$ with those fixed signs.

\begin{lemma}\label{lem:feasible-cubes}
Let $P,Q\in\cF$, $J\subseteq E_0$, and $S=\Sep(P,Q)\cup J$.
If $P^K\in\cF$ for every $K\subseteq S$, then setting the entries of
$P$ on $S$ to zero gives a feasible covector. In particular,
$J\in\Delta_{PQ}$.
If $(\M,0)$ is normal crossing, $J\in\Delta_{PQ}$ implies that every
$P^K$, $K\subseteq S$, is feasible.
\end{lemma}
\begin{proof}
The first assertion follows by strong elimination and induction on $|S|$.

Conversely, $J\in\Delta_{PQ}$ gives a feasible $G\leq P,Q$ with
$J\subseteq z(G)$, and hence $S\subseteq z(G)$. If $(\M,0)$ is normal
crossing, the restriction to $z(G)$ is free, so every sign assignment there
extends to a covector $Y$. The covector $G\circ Y$ fills those zero
entries and preserves all other signs, including the sign at $0$.
Choosing the signs on $z(G)\setminus S$ to agree with $P$ gives every
vertex $P^K$ of the cube.
\end{proof}

\begin{corollary}\label{cor:face-subsets}
Assume that $(\M,0)$ is normal crossing. The zero-set map is an order-reversing
bijection between the feasible faces $G\leq P$ and the simplices of
$\Delta_{PP}$. The face indexed by $J$ has zero entries on $J$ and
the signs of $P$ elsewhere. If $P,Q$ have a nonempty common face,
the covector of $\bar P\cap\bar Q$ has zero set $\Sep(P,Q)$.
\end{corollary}
\begin{proof}
Apply Lemma~\ref{lem:feasible-cubes} with $Q=P$ for the first
assertion. A face below $P$ is uniquely determined by its zero set.
For the second assertion, apply the lemma with $J=\varnothing$.
The resulting covector is the largest common face of $P,Q$.
\end{proof}

Assume that $(\M,0)$ is normal crossing. For $P,R\in\cF$, put
\[
 \begin{aligned}
 L_P^R&=\{z(G)\mid G\leq P\text{ feasible},\ z(G)\subseteq\Sep(P,R)\}\\
      &=\{J\in\Delta_{PP}\mid J\subseteq\Sep(P,R)\}.
 \end{aligned}
\]
The equality follows from Corollary~\ref{cor:face-subsets}, so this
is a simplicial complex.
Its vertices label the facets of $\bar P$ separating $P$ from $R$.
For hyperplanes these are the facets visible from a point of $R$.
Their union projects onto a convex set, and their nonempty
intersections are convex. The nerve lemma therefore gives
contractibility when $P\ne R$.

\begin{lemma}\label{lem:separating-facets}
If $(\M,0)$ is normal crossing, then $L_P^R$ is contractible for $P\ne R$.
\end{lemma}
\begin{proof}
The face-filter lemma \cite[Lemma 4.3.11]{OrientedMatroids1999}
says that for topes $P\ne\pm R$ the poset
\[
 \{G\mid 0<G<P,\ z(G)\subseteq\Sep(P,R)\}
\]
has contractible order complex. Here the order complex has a vertex
for each covector and a simplex for each chain of covectors.
Since $P(0)=R(0)=+$, these covectors are feasible and $P\ne\pm R$.
The zero-set map identifies this poset with the nonempty simplices of
$L_P^R$, ordered by reverse inclusion. Its order complex is therefore
the barycentric subdivision of $L_P^R$, which has the same underlying
space. The assertion follows.
\end{proof}

The last lemma converts a weighted count of hyperplane intersections
meeting a face into a count of its bounded faces. For a ray in
dimension one, it reads $(1-t)+t=1$, with only the endpoint contributing
on the right.

\begin{lemma}\label{lem:face-enumeration}
If $(\mathcal M,0)$ is normal crossing, for $P,Q\in\cF$ and $F=\bar P\cap\bar Q$,
\[
 \sum_{J\in\Delta_{PQ}}t^{|J|}(1-t)^{r-|J|}=h(F_\bd,t).
\]
\end{lemma}
\begin{proof}
We may assume $F\ne\varnothing$.
Expanding $t^{|J|}=(1+(t-1))^{|J|}$ gives
\[
 \sum_{J\in\Delta_{PQ}}t^{|J|}(1-t)^{r-|J|}
 =\sum_{I\in\Delta_{PQ}}c_I(t-1)^{r-|I|},
\text{ where }
 c_I=\sum_{\substack{J\in\Delta_{PQ}\\I\subseteq J}}(-1)^{r-|J|}.
\]

If $z(F)\not\subseteq I$, choose $i\in z(F)\setminus I$. Since
$F\subseteq H_i$, we have $
 F\cap H_J=F\cap H_{J\cup\{i\}}.$
Hence the terms in $c_I$ cancel in pairs indexed by $J$
and $J\cup\{i\}$, with $i\notin J$, giving $c_I=0$.

If $z(F)\subseteq I$, Corollary~\ref{cor:face-subsets} gives a unique
feasible face $G\leq F$ with $z(G)=I$. The sets $J\in\Delta_{PQ}$
containing $I$ correspond to the feasible faces $H\leq G$.
For these faces, normal crossing gives $\dim H=r-|z(H)|$, so
\[
 c_I=\sum_{\substack{H\leq G\\H \text{ feasible}}}(-1)^{\dim H}.
\]
Denote by $G^\infty\leq G$ the largest covector below $G$ with
$G^\infty(0)=0$, obtained by composing all such covectors. Then
$H\leq G$ is not feasible if and only if $H\leq G^\infty$.
Moreover, $G$ is bounded if and only if $G^\infty=\mathbf0$.

We will make use of the Euler relation for covectors, see 
\cite[Corollary 4.3.8]{OrientedMatroids1999}, which states that for every nonzero covector $G'$
\[
 \sum_{\mathbf0<H\leq G'}(-1)^{\dim H}=1.
\]
Applying this to $G'=G$ and splitting the sum, we obtain
\[
 \begin{aligned}
  1 &=\sum_{\mathbf0<H\leq G}(-1)^{\dim H}  = c_I
    +\sum_{\mathbf0<H\leq G^\infty }
       (-1)^{\dim H}.
 \end{aligned}
\]
If $G$ is bounded $G^\infty=\mathbf{0}$ and the second summand is zero, so $c_I=1$. Otherwise, the Euler identity for $G'=G^\infty$ applies to the second summand so that $c_I=0$.

Consequently, we obtain
\[
 \sum_{J\in\Delta_{PQ}}t^{|J|}(1-t)^{r-|J|}
 =\sum_{\substack{G\leq F\\G\text{ bounded}}}(t-1)^{\dim G}
 =h(F_\bd,t)\qedhere
\]
\end{proof}

\section{Koszul duality and the inverse Varchenko matrix}\label{sec:koszul}

Let $(\M,0)$ be an affine oriented matroid. We use the definitions
\[
 V_b(P,Q)=b_{\Sep(P,Q)},\qquad
 W_q(P,Q)=(-q)^{d(P,Q)}h((\bar P\cap\bar Q)_\bd,q^2)
\]
for all $P,Q\in\cF$, interpreting faces as in
Section~\ref{sec:background}. We denote the submatrix of $W_q$ indexed by bounded topes by $S_q$. The variables $b_i$ in this section are indexed by $E_0$. The main goal of this section is to relate the matrices $V_b$ and $W_q$ in the case that $(\M,0)$ is normal crossing via categorification and Koszul duality.

\subsection{Taut path and intersection algebras}
Fix a field $k$. We compose paths from left to right, so $e_PAe_Q$
has initial region $P$ and terminal region $Q$.

\begin{definition}
The \emph{taut path algebra} is the vector space
\[
 \Taut=\bigoplus_{P,Q\in\cF}k\tau_{PQ}
\]
with multiplication
\[
 \tau_{PQ}\tau_{SR}=
 \begin{cases}
  \tau_{PR},&Q=S\text{ and }d(P,Q)+d(Q,R)=d(P,R),\\
  0,&\text{otherwise}.
 \end{cases}
\]
\end{definition}

For each $P$, the element $e_P=\tau_{PP}$ is idempotent, and
$\sum_{P\in\cF}e_P=1$. A \emph{path} is a sequence of regions
sharing successive facets. It is \emph{taut} if it crosses each
separating hyperplane once and no others. The element $\tau_{PQ}$
is represented by any taut path from $P$ to $Q$.

\begin{definition}
For $P,Q\in\cF$, let
\[
 R_{PQ}=k[x_i\mid i\in E_0]/(x_J\mid J\notin\Delta_{PQ}),
 \qquad x_J=\prod_{i\in J}x_i.
\]
Write $1_{PQ}$ for the class of $1$ in this block. We use
$x_\varnothing=1$, so $R_{PQ}=0$ when $P,Q$ have no common feasible
face. The \emph{chamber intersection algebra} is
\[
 \Inter=\bigoplus_{P,Q\in\cF}R_{PQ},
\]
where the variables act centrally, products with different middle
labels are zero, and
\begin{equation}\label{eq:intersection-product}
 1_{PQ}1_{QR}=x_{\Sep(P,Q)\cap\Sep(Q,R)}1_{PR}.
\end{equation}
\end{definition}
The products respect the defining ideals, since a common face of $P,R$
on the hyperplanes in $\Sep(P,Q)\cap\Sep(Q,R)$ is also a face of
$Q$. Moreover, one may show associativity by counting repeated crossings in each coordinate.

Again $e_P=1_{PP}$ is idempotent and $\sum_Pe_P=1$, so the two
algebras have the same degree-zero algebra $\bigoplus_Pke_P$.
For adjacent regions separated by $H_i$,
\[
 1_{PQ}1_{QP}=x_ie_P,\qquad \tau_{PQ}\tau_{QP}=0.
\]
\begin{remark}
For realizable normal crossing arrangements, $\Inter$ is the algebra
$\widetilde B_{\mathrm{ext}}$ of \cite[Section 4.5]{BLPWgale}.
Their $\widetilde B$ is its corner on the feasible chambers bounded
with respect to a chosen generic objective function. Imposing the
linear parameter relations then gives their $B$
\cite[Definition 4.3 and Remark 4.5]{BLPWgale}.
Kowalenko--Mautner use the same face rings and multiplication for
oriented matroid programs \cite[Definition 5.9]{KowalenkoMautner}.
Our algebra uses all feasible chambers and retains the parameters.
Boundedness for an objective function differs from the boundedness
defining $\cB$ here.
\end{remark}

Both algebras are graded by $\NN^{E_0}$, with
\[
 \mdeg\tau_{PQ}=\mdeg1_{PQ}=\one_{\Sep(P,Q)},\qquad
 \mdeg x_i=2\mathbf e_i.
\]
Here $\one_J$ is the indicator vector of $J$, and $\mathbf e_i$ is
the $i$th standard basis vector. Total degree is the sum of the
coordinates of the multidegree.

\subsection{Bases and Hilbert matrices}
The block $e_P\Taut e_Q$ has basis $\{\tau_{PQ}\}$, and
$e_P\Inter e_Q$ has basis
\[
 \{x^\nu1_{PQ}\mid\nu\in\NN^{E_0},\ \supp\nu\in\Delta_{PQ}\}
\]
since the defining monomial ideal of $R_{PQ}$ is spanned by precisely
the monomials whose supports do not belong to $\Delta_{PQ}$.
For either algebra $A$, define its multigraded Hilbert matrix by
\[
 \Hilb_b(A)_{P,Q}=\sum_{\mathbf m\in\NN^{E_0}}
       \dim_k(e_PA_{\mathbf m}e_Q)b^{\mathbf m},
 \qquad b^{\mathbf m}=\prod_{i\in E_0}b_i^{m_i}.
\]
The bases give
\begin{align*}
 \Hilb_b(\Taut)_{P,Q}&=b_{\Sep(P,Q)}=V_b(P,Q),\\
 \Hilb_b(\Inter)_{P,Q}&=b_{\Sep(P,Q)}
    \sum_{J\in\Delta_{PQ}}\prod_{i\in J}\frac{b_i^2}{1-b_i^2}.
\end{align*}
By applying Lemma~\ref{lem:face-enumeration} to this formula, we obtain the following relation with the $q$-intersection matrix.
\begin{corollary}\label{cor:graded-hilbert}
For total degree, $\Hilb_q(\Taut)=V_q$. If $(\M,0)$ is normal crossing, then
\[
 \Hilb_q(\Inter)=(1-q^2)^{-r}W_{-q}.
\]
\end{corollary}
\subsection{Inverse Varchenko matrix via Koszul duality}
Let $A$ be a nonnegatively graded $k$-algebra with finite-dimensional
graded pieces and $A_0=\bigoplus_{i\in I}ke_i$ for pairwise orthogonal
idempotents $e_i$. It is \emph{quadratic} if
$A=T_{A_0}(V)/(\mathcal R)$, where $V=A_1$ and
$\mathcal R\subseteq V\otimes_{A_0}V$. Its quadratic dual is
\[
 A^!=T_{A_0}(V^*)/(\mathcal R^\perp),
\]
where the bimodule dual reverses endpoints. Dual generators retain
their multidegrees \cite[Definition 2.8.1]{BeilinsonGinzburgSoergel1996}.
The algebra is \emph{Koszul} if $A_0$ has a graded projective resolution
whose term in homological degree $n$ is generated in total degree $n$.
For a Koszul algebra the quadratic dual is also Koszul, and
\[
 \Hilb_b(A)\Hilb_{-b}(A^!)^{\mathsf T}=\mathrm{Id}_I
\]
by \cite[Proposition 2.9.1 and Lemma 2.11.1]{BeilinsonGinzburgSoergel1996}.

\begin{theorem}\label{thm:chamber-koszul}
For a normal crossing affine oriented matroid, $\Taut$ and $\Inter$
are quadratic and Koszul over every field, and $\Taut^!\cong\Inter$.
\end{theorem}

\begin{remark}
Kowalenko--Mautner ask whether their category $\mathcal O$ algebras
are Koszul for every generic oriented matroid program
\cite[Section 1.3]{KowalenkoMautner}. They prove Koszulity in the
Euclidean case and the numerical Koszul identity in general
\cite[Theorems 7.13--7.14]{KowalenkoMautner}.
Theorem~\ref{thm:chamber-koszul} requires no Euclidean hypothesis,
but concerns the algebras on all feasible chambers. Passing to the
category $\mathcal O$ algebras involves a further choice of chambers
and linear parameter relations, so their Koszulity question requires
an additional argument.
\end{remark}

\begin{corollary}\label{cor:multivariate-inverse}
If $(\M,0)$ is normal crossing, then $V_b^{-1}=\Hilb_{-b}(\Inter)$. Explicitly,
\[
 (V_b^{-1})(P,Q)=(-1)^{d(P,Q)}b_{\Sep(P,Q)}
     \sum_{J\in\Delta_{PQ}}\prod_{i\in J}\frac{b_i^2}{1-b_i^2}.
\]
\end{corollary}

\begin{corollary}\label{thm:varinv}
For a normal crossing affine oriented matroid,
\[
 W_q=(1-q^2)^rV_q^{-1},\qquad
 S_q=(1-q^2)^r(V_q^{-1})^{\cB}.
\]
\end{corollary}

\begin{example}\label{ex:two-points}
For the two-point arrangement of Example~\ref{ex:line},
with its three chambers in order,
\[
 V_b=\begin{pmatrix}
  1&b_1&b_1b_2\\ b_1&1&b_2\\ b_1b_2&b_2&1
 \end{pmatrix},\qquad
 (V_b^{-1})^{\cB}
 =\frac{1-b_1^2b_2^2}{(1-b_1^2)(1-b_2^2)}.
\]
The unique bounded chamber is an interval. Setting $b_1=b_2=q$ gives
\[
 S_q=1+q^2,\qquad
 (V_q^{-1})^{\cB}=\frac{1+q^2}{1-q^2}.
\]
\end{example}

Combining Proposition~\ref{prop:Betti-Varchenko} with
Corollaries~\ref{cor:multivariate-inverse} and~\ref{thm:varinv} gives
an explicit formula for the Betti matrix at $b_0=0$ when $(\M,0)$
is normal crossing.

\begin{corollary}\label{cor:Betti-entries}
If $(\M,0)$ is normal crossing, then for all $P,Q\in\cF$,
\[
 \left.B_b(P,Q)\right|_{b_0=0}
 =(-1)^{r+d(P,Q)}b_{\Sep(P,Q)}
   \sum_{J\in\Delta_{PQ}}\prod_{i\in J}\frac{b_i^2}{1-b_i^2}.
\]
In particular,
\[
 \left.B_b\right|_{b_0=0,\ b_i=q\ (i\in E_0)}
 =(-1)^r(1-q^2)^{-r}W_q.
\]
\end{corollary}

\begin{remark}
If $(\M,0)$ is normal crossing, the Betti matrix at $b_0=0$ is recovered from
the multigraded Hilbert matrix of the intersection algebra $\Inter$:
\[
 \left.B_b\right|_{b_0=0}=(-1)^r\Hilb_{-b}(\Inter).
\]
It would be interesting to extend this categorification beyond normal
crossing, replacing $\Inter$ by the dg Koszul dual of $\Taut$.
We also expect an extension to general $b_0$, giving a categorification
of Lam's full Betti matrix.
\end{remark}

\begin{remark}\label{ex:concurrent-lines}
Consider three distinct lines in the plane meeting at a point.
Their chamber graph is a hexagon. The two taut paths between
opposite chambers have length three and represent the same element
of $\Taut$. The only quadratic relations kill immediate backtracks,
so they cannot identify these two paths. Thus $\Taut$ is not quadratic.

Moreover, for each chamber $P$,
\[
 (W_qV_q)(P,P)=1-2q^2+2q^4-q^6\ne(1-q^2)^2,
\]
so the normal crossing hypothesis in Corollary~\ref{thm:varinv}
cannot be omitted. Lam's face-flag formula still computes $V_b^{-1}$
through his Betti matrix, as recalled in
Proposition~\ref{prop:Betti-Varchenko}.
\end{remark}

\subsection{Proof of quadratic duality}
Let
\[
 \Cube_1=k\bigl(+\mathop{\rightleftarrows}^{u}_{v}-\bigr),\qquad
 \theta=uv+vu,\qquad
 \CubeZero_1=\Cube_1/(\theta)=\Cube_1/(uv,vu).
\]
Identifying the
reversed dual arrows with the arrows having the same endpoints gives
\[
 \Cube_1^!\cong\CubeZero_1,\qquad
 \CubeZero_1^!\cong\Cube_1.
\]
Take one copy for each hyperplane and consider the algebras
\[
 \Cube=\bigotimes_{i\in E_0}\Cube_1,\qquad
 \CubeZero=\bigotimes_{i\in E_0}\CubeZero_1
          =\Cube/(\theta_i\mid i\in E_0),
\]
where $\theta_i$ is the central element $\theta$ in the $i$th factor.
Quadratic duality turns commuting squares into anticommuting
squares, so the dual of an ordinary tensor product is initially a
signed tensor product. Rescaling the arrows removes these signs.
Order $E_0$, write $\sigma^i$ for $\sigma$ with its $i$th sign reversed,
and put
\begin{equation}\label{eq:cube-signs}
 \epsilon_i(\sigma)=(-1)^{|\{j<i\mid\sigma(j)=-\}|}.
\end{equation}
For $i\ne j$,
$\epsilon_i(\sigma)\epsilon_j(\sigma^i)
 =-\epsilon_j(\sigma)\epsilon_i(\sigma^j)$.
Thus mapping the dual arrow $\sigma\to\sigma^i$ to
$\epsilon_i(\sigma)$ times the corresponding arrow of $\Cube$
changes anticommuting squares into commuting squares and gives
\[
 \CubeZero^!\cong\Cube.
\]

\begin{proposition}\label{prop:feasible-models}
Set $e=\sum_{P\in\cF}e_P$. Then
\[
 \Taut\cong e\CubeZero e.
\]
If the arrangement, or affine oriented matroid, is normal crossing,
then also
\[
 \Inter\cong\Cube/\Cube(1-e)\Cube.
\]
\end{proposition}
\begin{proof}
The first identification is left to the reader.
For the second, we compare bases in each block.
Fix $P,Q\in\cF$, and let $\widehat p_{PQ}$ be the taut path class
in $\Cube$. The block $e_P\Cube e_Q$ has basis
\[
 \bigl\{\theta^\nu\widehat p_{PQ}\bigm|\nu\in\NN^{E_0}\bigr\}.
\]
The quotient by $\Cube(1-e)\Cube$ keeps exactly those basis
monomials for which the cube obtained from $P$ by varying the signs
in $\Sep(P,Q)\cup\operatorname{supp}\nu$ has only feasible vertices.

If $(\M,0)$ is normal crossing, Lemma~\ref{lem:feasible-cubes}, applied with
$J=\operatorname{supp}\nu$, identifies this condition with
\[
 \bar P\cap\bar Q\cap H_{\operatorname{supp}\nu}\ne\varnothing,
 \qquad\text{equivalently,}\qquad
 \operatorname{supp}\nu\in\Delta_{PQ}.
\]
The surviving monomials have the same exponent set as the basis
of $R_{PQ}$. Hence the map
$\theta^\nu\widehat p_{PQ}\longmapsto x^\nu 1_{PQ}$ yields the second claimed isomorphism.
\end{proof}

Assume now that $(\M,0)$ is normal crossing. It remains to show that
$e\CubeZero e$ is quadratic.
Every pair of feasible topes is joined by a taut path, so the
arrows between adjacent regions generate $\Taut$, hence the algebra is generated in degree one. To see that quadratic relations suffice, note that any two taut paths with the same endpoints are related by
feasible square moves \cite[Proposition 3.5]{KowalenkoMautner}.
For a non-taut path, square moves expose an immediate backtrack
\cite[proof of Proposition 3.6]{KowalenkoMautner}. These path arguments
use only the feasible topes and the normal crossing condition. Thus commuting
feasible squares and killing immediate backtracks give a quadratic
presentation of $\Taut$.
The corner--quotient lemma \cite[Lemma 8.22]{BradenLicataProudfootWebster2012}
now gives
\[
 \Taut^!
 \cong\CubeZero^!/\CubeZero^!(1-e)\CubeZero^!
 \cong\Cube/\Cube(1-e)\Cube
 \cong\Inter,
\]
proving that the two algebras are quadratic duals.
\subsection{Proof of Koszulity}
Assume that $(\M,0)$ is normal crossing. We adapt the method of the corrected proof
of \cite[Lemma 8.25]{BradenLicataProudfootWebster2012} to all feasible
chambers. We decompose the Koszul complex into multidegrees and use
Lemma~\ref{lem:separating-facets} for the required acyclicity.
Since $\Inter^!\cong\Taut$, it suffices to prove
exactness of the augmented Koszul complex of $\Inter$
\cite[Theorem 2.6.1, \S2.8]{BeilinsonGinzburgSoergel1996}.
Fix $P,R\in\cF$. Multiplying the Koszul complex by the corresponding
idempotents gives
\[
 K_n=e_P\Inter\otimes_{\Inter_0}(\Taut_n)^*e_R
 =\bigoplus_{\substack{Q\in\cF\\d(Q,R)=n}}
     e_P\Inter e_Q\otimes_k k\tau_{RQ}^*.
\]
The differential preserves multidegree, so we can further fix
$\mathbf m\in\NN^{E_0}$ and work in the subcomplex
$K_{\bullet,\mathbf m}$.

To prove exactness, we distinguish two cases. If some $m_i\geq2$,
one component of the differential pairs the basis vectors.
If all $m_i\leq1$, we identify the nonzero component with the
reduced cochain complex of $L_P^R$, defined before
Lemma~\ref{lem:separating-facets}.

Using the signs in \eqref{eq:cube-signs}, we decompose the differential as
$\partial=\sum_{i\in E_0}\partial_i$, where
\[
 \partial_i(z\otimes\tau_{RQ}^*)=
 \begin{cases}
  \epsilon_i(Q)\,z 1_{Q,Q^i}\otimes\tau_{R,Q^i}^*,
     &i\in\Sep(Q,R),\ Q^i\in\cF,\\
  0,&\text{otherwise}.
 \end{cases}
\]
Moreover, we will denote
\begin{equation}\label{eq:koszul-exponent}
 \nu(Q)=\tfrac12\bigl(\mathbf m-\one_{\Sep(P,Q)}
                              -\one_{\Sep(Q,R)}\bigr),
\end{equation}
so that $K_{n,\mathbf m}$ has basis vectors
\[
 v_Q=x^{\nu(Q)}1_{PQ}\otimes\tau_{RQ}^*,\qquad
 \begin{cases}
  Q\in\cF,\quad d(Q,R)=n,\\
  \nu(Q)\in\NN^{E_0},\\
  \operatorname{supp}\nu(Q)\in\Delta_{PQ}.
 \end{cases}
\]
The two indicator entries in \eqref{eq:koszul-exponent} sum to
$1$ where $P(i)\ne R(i)$, and to $0$ or $2$ otherwise.
The component $K_{n,\mathbf m}$ is zero unless
$\mathbf m\equiv\one_{\Sep(P,R)}\pmod2$. Hence we assume this parity condition from now on.
For $i\in\Sep(Q,R)$, \eqref{eq:koszul-exponent} and
\eqref{eq:intersection-product} give
\[
 \nu(Q^i)-\nu(Q)=
 \begin{cases}
  \mathbf e_i,&P(i)=R(i),\\
  0,&P(i)\ne R(i),
 \end{cases}
 \qquad
 \partial_i v_Q=
 \begin{cases}
  \epsilon_i(Q)v_{Q^i},&\text{if }v_{Q^i}\text{ exists},\\
  0,&\text{otherwise}.
 \end{cases}
\]

\smallskip
\noindent\emph{Case 1: $m_i\geq2$ for some $i$.}
The basis vectors occur in pairs $v_Q,v_{Q^i}$, with
$\partial_i v_Q=\epsilon_i(Q)v_{Q^i}$ when $Q(i)\ne R(i)$.
Since $m_i\geq2$, formula~\eqref{eq:koszul-exponent} gives
$\nu(Q^i)\in\NN^{E_0}$ for every basis vector $v_Q$, and
\[
 \Sep(P,Q)\cup\operatorname{supp}\nu(Q)
 =\Sep(P,Q^i)\cup\operatorname{supp}\nu(Q^i).
\]
By Lemma~\ref{lem:feasible-cubes}, the cube based at $P$ with this
common set of directions is feasible. It contains $Q^i$, and the
same lemma gives $\operatorname{supp}\nu(Q^i)\in\Delta_{P,Q^i}$.
Hence, every partner $v_{Q^i}$ exists.
In each homological degree, set
\[
 U_n=\langle v_Q\mid Q(i)\ne R(i)\rangle_k,\qquad
 V_n=\langle v_Q\mid Q(i)=R(i)\rangle_k.
\]
The pairing gives $\partial_i:U_n\xrightarrow{\sim}V_{n-1}$.
Every $\partial_j$ with $j\ne i$ preserves $Q(i)$, hence $U$ and $V$.
For a cycle $z=u+v\in U_n\oplus V_n$, choose
$y\in U_{n+1}$ with $\partial_i y=v$.
Then $w=z-\partial y$ is a cycle in $U_n$.
The $V_{n-1}$-component of $\partial w$ is $\partial_i w$, so
injectivity of $\partial_i$ gives $w=0$ and $z=\partial y$.

\smallskip
\noindent\emph{Case 2: every $m_i\leq1$.}
The parity condition gives
$\mathbf m=\one_{\Sep(P,R)}$.
A basis vector $v_Q$ now has $\nu(Q)=0$ and
\[
 \Sep(P,Q)\uplus\Sep(Q,R)=\Sep(P,R).
\]
Let $G$ be the covector of the nonempty
common face $\bar P\cap\bar Q$. Corollary~\ref{cor:face-subsets}
gives $z(G)=\Sep(P,Q)$.
Associate $v_Q$ with the simplex $z(G)$ of $L_P^R$.
The inverse recovers $Q=G\circ R$, a feasible tope with
$G\leq P,Q$ and $\Sep(P,Q)=z(G)$.
The basis of $K_{n,\mathbf m}$ therefore corresponds to the simplices
with $|z(G)|=d(P,R)-n$.

A nonzero differential term in direction $i$ sets $G(i)$ to zero,
so it adds $i$ to the simplex $z(G)$.
Rescaling $v_Q$ by $\prod_{j\in z(G)}\epsilon_j(P)$
gives the usual simplicial coboundary signs. This component is the reduced
cochain complex
\[
 K_{n,\mathbf m}
 \cong\widetilde C^{\,d(P,R)-n-1}(L_P^R;k).
\]
For $P\ne R$, Lemma~\ref{lem:separating-facets} makes it exact.
For $P=R$, only $v_P$ remains, and the augmentation
maps it isomorphically to $k$.

\section{The \texorpdfstring{$q$}{q}-intersection determinant}\label{sec:qdet}
\begin{proof}[Proof of Theorem~\ref{thm:intersection-real}]
Let $\C$ be an essential real affine hyperplane arrangement with normal crossings.  By Theorem~\ref{thm:Lam-real},
\[
 \det\bigl((V_b^{-1})^{\cB}\bigr)
 =\frac{\displaystyle\prod_{F\in\cL_0\setminus\{E\}}
          (1-b_{E\setminus F}^2)^{w_F}}
        {\displaystyle\prod_{F\in\cL\setminus(\cL_0\cup\{\varnothing\})}
          (1-b_F^2)^{w_F}}.
\]

By Proposition~\ref{prop:normal-crossing-ranks}
and Lemma~\ref{lem:free-coextension},
the flats avoiding $0$ contribute only when they are singletons,
and the proper flats containing $0$ are $F=K\cup\{0\}$ for
$K\in\cL(M_0)\setminus\{E_0\}$. Put
\[
 \lambda_K=w_{K\cup\{0\}}
       =\beta(M_{0,K})\mup((M_0^K)^*).
\]
Specializing $b_i$ to $q$, for all $i\in E_0$, in the expression above gives
\[
 \det\bigl((V_q^{-1})^{\cB}\bigr)
 =\frac{\displaystyle\prod_{K\in\cL(M_0)\setminus\{E_0\}}
          (1-q^{2|E_0\setminus K|})^{\lambda_K}}
        {(1-q^2)^{\sum_{i\in E_0}w_{\{i\}}}}.
\]
The inverse-matrix identity of Corollary~\ref{thm:varinv} and
$|\cB|=\beta(M)$ therefore give
\[
 \det S_q
 =(1-q^2)^{r\beta(M)-\sum_{i\in E_0}w_{\{i\}}}
   \prod_{K\in\cL(M_0)\setminus\{E_0\}}
      (1-q^{2|E_0\setminus K|})^{\lambda_K}.
\]
Corollary~\ref{cor:degree} identifies the exponent outside the product
with $-\sum_K\lambda_K$. Hence
\[
 \det S_q
 =\prod_{K\in\cL(M_0)\setminus\{E_0\}}
    \left(\frac{1-q^{2|E_0\setminus K|}}{1-q^2}\right)^{\lambda_K}.
\]
If $K$ is not cyclic, $M_0^K$ has a coloop, so its dual has a loop
and $\mup((M_0^K)^*)=0$. The product reduces to
\eqref{eq:Sq-determinant}.
\end{proof}

The same calculation, using the oriented-matroid form of
Corollary~\ref{thm:varinv} and Lemma~\ref{lem:free-coextension},
proves the following reduction independently of realizability.

\begin{corollary}\label{cor:conjecture-reduction}
For a normal crossing affine oriented matroid, the conjectural determinant
formula \eqref{eq:Lam-determinant} implies \eqref{eq:Sq-determinant}.
\end{corollary}

\begin{example}\label{ex:four-lines}
Consider the two normal crossing arrangements $\C$ and $\C'$ below,
each consisting of four lines with $H_1\parallel H_2$. Put
$E_0=\{1,2,3,4\}$ and order their bounded chambers as
$\cB=(P_1,P_2)$ and $\cB'=(P'_1,P'_2)$, respectively.

\begin{center}
 \begin{minipage}{0.46\textwidth}
  \centering
  \begin{tikzpicture}[scale=0.8]
   \fill[gray!15] (0,0) -- (0,1) -- (1,1) -- cycle;
   \fill[gray!15] (0,0) -- (0,-1) -- (-1,-1) -- cycle;
   \draw[thick] (-2,1) -- (2,1) node[right] {$H_1$};
   \draw[thick] (-2,-1) -- (2,-1) node[right] {$H_2$};
   \draw[thick] (0,-2) -- (0,2) node[above] {$H_3$};
   \draw[thick] (-2,-2) -- (2,2) node[above right] {$H_4$};
   \node at (0.3,0.7) {$P_1$};
   \node at (-0.3,-0.7) {$P_2$};
  \end{tikzpicture}
  \par\smallskip
  $\C$
 \end{minipage}\hfill
 \begin{minipage}{0.46\textwidth}
  \centering
  \begin{tikzpicture}[scale=0.8]
   \fill[gray!15] (0,0) -- (0,-1) -- (-1,-1) -- cycle;
   \fill[gray!15] (0,-1) -- (0,-2) -- (-2,-2) -- (-1,-1) -- cycle;
   \draw[thick] (-3,-2) -- (1.5,-2) node[right] {$H_1$};
   \draw[thick] (-3,-1) -- (1.5,-1) node[right] {$H_2$};
   \draw[thick] (0,-3) -- (0,1) node[above] {$H_3$};
   \draw[thick] (-3,-3) -- (1,1) node[above right] {$H_4$};
   \node at (-0.3,-0.7) {$P'_1$};
   \node at (-0.7,-1.5) {$P'_2$};
  \end{tikzpicture}
  \par\smallskip
  $\C'$
 \end{minipage}
\end{center}

In $\C$, the two bounded chambers are triangles whose closures meet
at a vertex. In $\C'$, they are a triangle and a quadrilateral whose
closures share an edge. Corollary~\ref{cor:multivariate-inverse} gives
\[
 \begin{aligned}
 \bigl(V_b(\C)^{-1}\bigr)^{\cB}
 ={}&(1-b_3^2)^{-1}(1-b_4^2)^{-1}\\
 &{}\cdot
 \begin{pmatrix}
  (1-b_{\{1,3,4\}}^2)(1-b_1^2)^{-1}&b_3b_4\\
  b_3b_4&(1-b_{\{2,3,4\}}^2)(1-b_2^2)^{-1}
 \end{pmatrix}
 \end{aligned}
\]
and
\[
 \begin{aligned}
 \bigl(V_b(\C')^{-1}\bigr)^{\cB'}
 ={}&(1-b_2^2)^{-1}(1-b_3^2)^{-1}(1-b_4^2)^{-1}\\
 &{}\cdot
 \begin{pmatrix}
  1-b_{\{2,3,4\}}^2&-b_2(1-b_{\{3,4\}}^2)\\
  -b_2(1-b_{\{3,4\}}^2)&
  (1-b_{\{1,2\}}^2)(1-b_{\{3,4\}}^2)(1-b_1^2)^{-1}
 \end{pmatrix}.
 \end{aligned}
\]
Both determinants are
\[
 \frac{(1-b_{E_0}^2)(1-b_{\{3,4\}}^2)}
 {(1-b_1^2)(1-b_2^2)(1-b_3^2)^2(1-b_4^2)^2}.
\]
This agrees with Theorem~\ref{thm:Lam-real}. For both arrangements,
the nonzero weights of proper nonempty flats of $M$ are
\[
 w_{\{0\}}=w_{\{0,1,2\}}=w_{\{1\}}=w_{\{2\}}=1,
 \qquad w_{\{3\}}=w_{\{4\}}=2.
\]

Setting $b_i=q$ and multiplying by $(1-q^2)^2$ gives
\[
 S_q(\C)=
 \begin{pmatrix}
  1+q^2+q^4&q^2\\
  q^2&1+q^2+q^4
 \end{pmatrix},\qquad
 S_q(\C')=
 \begin{pmatrix}
  1+q^2+q^4&-q(1+q^2)\\
  -q(1+q^2)&(1+q^2)^2
 \end{pmatrix}.
\]
Thus
\[
 \det S_q(\C)=\det S_q(\C')
 =(1+q^2)^2(1+q^4)=\qint{2}\,\qint{4}.
\]
The two factors correspond to the proper cyclic flats
$\{1,2\}$ and $\varnothing$ of $M_0$, each with exponent one in
Theorem~\ref{thm:intersection-real}.
\end{example}

\section*{Acknowledgements}
We thank Thomas Lam, Manfred Lehn, and Duco van Straten for
interesting discussions.

We used LLMs for an extensive literature search on determinant formulas
and to assist with computations in many concrete examples.

\end{document}